\documentclass{article}
 \usepackage{amsmath}
 \usepackage{pifont}

\newtheorem{thm}{Theorem}[section]

\newtheorem{prop}[thm]{Proposition}

\newtheorem{lem}[thm]{Lemma}
\newtheorem{Def}[thm]{Definition}
\newtheorem{rem}[thm]{Remark}

\newtheorem{ex}[thm]{Example}

\newcommand{\be}{\begin{equation}}
\newcommand{\ee}{\end{equation}}
\newcommand{\ben}{\begin{enumerate}}
\newcommand{\een}{\end{enumerate}}
\newcommand{\beq}{\begin{eqnarray}}
\newcommand{\eeq}{\end{eqnarray}}
\newcommand{\beqn}{\begin{eqnarray*}}
\newcommand{\eeqn}{\end{eqnarray*}}

\newcommand{\pa}{\partial}

\newcommand{\qed}{\hspace*{\fill}Q.E.D.}  

\begin{document}
\title{On  Sprays of Scalar Curvature and Metrizability}
\author{Guojun Yang }
\date{}
\maketitle
\begin{abstract}

 Every Finsler metric naturally induces a spray but not so for the converse.
 The notion for sprays of scalar (resp. isotropic) curvature has been known as a generalization
 for
  Finsler metrics of scalar (resp. isotropic) flag curvature. In this
 paper, a new notion, sprays of constant curvature, is introduced and especially it  shows that a spray
 of isotropic curvature is not necessarily of constant curvature even in dimension $n\ge3$. Further,
  complete conditions are given   for sprays of isotropic (resp. constant) curvature to be
 Finsler-metrizabile.
 As applications of such a result, the local  structure
is determined  for locally projectively flat Berwald sprays of
isotropic (resp. constant) curvature which are Finsler-metrizable,
and some more
  sprays of isotropic
 curvature are  discussed for  their metrizability.
Besides, the metrizability problem is also investigated for
 sprays of scalar curvature under certain curvature conditions.

  \bigskip
\noindent {\bf Keywords:}  Finsler Metric, Spray, Berwald Spray,
Metrizability, Scalar/Isotropic/Constant Curvature, Projective
Flatness

\noindent
 {\bf 2010 Mathematics Subject Classification: }
53C60,  53B40
\end{abstract}

\section{Introduction}

Spray geometry  studies the properties of  spray manifolds, a kind
of  path spaces. Spray geometry is more general than Finsler
geometry, because every Finsler metric induces a natural spray but
there are a lot of sprays which cannot be induced by any Finsler
metric (\cite{EM, LMY, Yang1}).  A spray ${\bf G}$ on a manifold
$M$ is a family of compatible second order ODEs which define a
special vector filed on a conical region $\mathcal{C}$ of
 $TM\setminus \{0\}$ (an important case is $\mathcal{C}=TM\setminus \{0\}$). The integral curves
of ${\bf G}$ projected onto $M$ are called geodesics of {\bf G}.
Many basic curvatures, such as Riemann curvature, Ricci curvature,
Weyl curvature, Berwald curvature and Douglas curvature, appearing
in Finsler geometry, are actually defined in spray geometry via
the spray coefficients. Geodesics and these basic curvatures play
an important role in the study of spray geometry.

The metrizability problem for a spray ${\bf G}$ seeks for a
Finsler metric whose spray is just {\bf G}, or whose geodesics
coincide with that of ${\bf G}$. A weaker problem is to consider
the projective metrizability of a given spray {\bf G}, which aims
to look for a Finsler metric projectively related to {\bf G}. So a
natural question is   to determine whether a given spray is
(projectively) Finsler-metrizable or not under certain curvature
conditions.

In \cite{Ma1}, M. Matsumoto proves that any two-dimensional spray
is locally projectively Finsler-metrizable. More generally, any
spray of scalar curvature is locally projectively
Finsler-metrizable (\cite{BM, Cram}).
 In \cite{Mu}, Z. Muzsnay gives
some sprays which are not Finsler-metrizable under some conditions
satisfied by the holonomy distribution generated from the
horizontal vector fields of a spray. In \cite{Yang1}, the present
author constructs a class of sprays whose metrizable and
non-metrizable conditions are completely determined respectively.
Inspired by the sprays constructed in \cite{Yang1}, S.G. Elgendi
and Z. Muzsnay discuss a more general class of sprays and prove
the non-metrizability of such sprays by computing the dimension of
the holonomy distribution under certain conditions (\cite{EM}). In
\cite{LS}, B. Li and Z. Shen introduce the notion of sprays with
isotropic curvature and give some non-metrizability conditions for
locally projectively flat sprays, and shows that  a locally
projectively flat spray with vanishing Riemann curvature is
metrizable (cf. \cite{Shen5}). In \cite{LMY}, Y. Li, X. Mo and Y.
Yu give a class of locally projectively flat Berwald sprays which
are non-metrizable.

  In this paper, we are going to study some special properties and the
metrizability problem of some special classes of sprays: locally
projectively flat Berwald sprays,  sprays of scalar curvature
(resp. isotropic curvature, constant curvature). A locally
 projectively flat spray, which is always of scalar curvature, means that its geodesics are locally
 straight lines. A Berwald spray means that its spray coefficients
 $G^i=G^i(x,y)$ are quadratic in $y$. A
spray ${\bf G}$ is said to be of {\it scalar curvature} if  its
Riemann curvature $R^i_{\ k}$ satisfies
 \be\label{yg1}
  R^i_{\ k}=R\delta^i_k-\tau_ky^i,
  \ee
where $R=R(x,y)$ and $\tau_k=\tau_k(x,y)$ are some homogeneous
functions (\cite{Shen6}). If in (\ref{yg1}) there holds
$R_{.i}=2\tau_i$,
 then ${\bf G}$  is said to be of {\it isotropic curvature}
 (\cite{LS}). A spray ${\bf G}$ is said to be of
 {\it constant curvature}, a new notion we introduce in this
 paper, if ${\bf G}$ satisfies (\ref{yg1}) with
  $$
 \tau_{i;k}=0 \ (\ \Leftrightarrow \ R=\tau_k=0,\ or \ R_{;i}=0).
  $$
 In the above,  we
use $T_{i;j}$ and $T_{i.j}$ to denote respectively the horizontal
and vertical covariant derivatives of the tensor $T$ with respect
to Berwald connection of a given spray.
 Some basic properties for a spray of constant
curvature are given in Theorem \ref{th01} below, in which,
 it  especially shows that an $n(\ge 3)$-dimensional spray
 of isotropic curvature is not necessarily of constant curvature,
 which is different from the Finslerian case.

Now we show the following Theorems \ref{th1} and \ref{th001} on
the metrizability for sprays of isotropic and constant curvature
respectively.

\begin{thm}\label{th1}
 Let {\bf G} be an $n$-dimensional spray  of isotropic  curvature
 $R^i_{\ k}=R\delta^i_k-\frac{1}{2}R_{.k}$.
  \ben
    \item[{\rm (i)}] If $R=0$, then {\bf G}  is locally Finsler-metrizable.
    \item[{\rm (ii)}] If $R\ne 0$ and $R$ is not a Finsler metric,
    then {\bf G} is not Finsler-metrizable.
    \item[{\rm (iii)}]  If $R$ is  a Finsler metric, then
    {\bf G}  is (locally) Finsler-metrizable if and only if $R_{;i}=R\omega_i$ for some closed 1-form $\omega=\omega_i(x)dx^i$.
 In this case, we have
    \ben
      \item[{\rm (iiia)}] if $n\ge 3$, then  $R_{;i}=0$ (or $\omega=0$).
      \item[{\rm (iiib)}] if $\omega=0$, then {\bf G}  is induced by  the Finlser metric
      $R$ with the flag curvature ${\bf K}=1$.
       \item[{\rm (iiib)}] if $\omega\ne 0$, then {\bf G}  is induced by  the Finlser metric
      $R/\lambda$ with the flag curvature ${\bf K}=\lambda$, where $\lambda\ne0$ is given by $\omega_i=(\ln|\lambda|)_{;i}$.
      This case happens only in dimension $n=2$.
     \een
  \een
\end{thm}

\begin{thm}\label{th001}
 Let {\bf G} be an $n$-dimensional spray of constant curvature and $Ric$ be the Ricci curvature of {\bf G}. Then {\bf G} is
(locally) Finsler-metrizable iff. $Ric=0$ or $Ric$ is a Finsler
metric. In
 this case, the Finsler metric $L$ inducing {\bf G} has vanishing flag curvature or
 $L$ is given by $L=Ric/(n-1)$ with the flag curvature $1$.
\end{thm}

By Theorems \ref{th1} and \ref{th001}, for a given spray {\bf G}
of isotropic or constant curvature, it is easy to check directly
whether {\bf G} is (locally) Finsler-metrizable or not. Meanwhile,
if it is metrizable, the corresponding Finsler metric is easily
obtained if
  $Ric\ne 0$. The metric in Theorem \ref{th1} can  be multiplied by a
suitable non-zero constant if we need the metric to be positive.
Theorem \ref{th1}(i) shows that a spray  with vanishing Riemann
curvature (not necessary to be locally projectively flat) is
locally metrizable (cf. \cite{Shen5, LS}).

As  applications of Theorems \ref{th1} and \ref{th001}, we
generalize a result in \cite{Yang1} (Theorem \ref{th410} below),
and make a check on some spays  whether they are
Finsler-metrizable or not (see Sec.\ref{sec7} below). On the other
hand, we can use Theorem \ref{th1} to obtain the  local structure
for locally projectively flat Berwald sprays of isotropic
curvature when they are Finsler-metrizable.

\begin{thm}\label{th2}
 Let {\bf G} be a projectively flat Berwald spray
of isotropic curvature on an open set $U\subset R^n$ with
$Ric\ne0$ (on a conical region $\mathcal{C}(U)$). Then {\bf G} is
Finsler metrizable (on $\mathcal{C}(U)$) if and only if {\bf G}
can be expressed as
 \be\label{ygj1}
 G^i=Py^i,\ \ \ \ \ \   P:=-\frac{1}{2}\Big[\ln|x'Ax+\langle
  B,x\rangle+C|\Big]_{x^k}y^k,
 \ee
where $A\ne 0$ is a constant symmetric matrix, $B$ is a constant
vector
  and $C$ is a constant number satisfying certain condition  such that the following function $L$ is a
  metric (defined on $\mathcal{C}(U)$),
 \be\label{ygj2}
  L:=\frac{4(x'Ax+\langle
  B,x\rangle+C)y'Ay-(2x'Ay+\langle
  B,y\rangle)^2}{4(x'Ax+\langle
  B,x\rangle+C)^2}.
  \ee
In this case, {\bf G} is induced by the
 metric $L=Ric/(n-1)$ of constant sectional curvature 1.
\end{thm}

In \cite{Shen5}, there is  a general description of the
construction for locally projectively flat Finsler metrics with
constant flag curvature 1.  As a special case, putting
$2A=(\delta_{ij}),B=0,C=1/2$ in Theorem \ref{th2}, we obtain
 \be\label{ygj3}
  L=\frac{(1+|x|^2)|y|^2-\langle x,y\rangle^2}{(1+|x|^2)^2}.
  \ee
Theorem \ref{th2} also implies that if we take in (\ref{ygj1}),
$P=-\big[\ln\sqrt{|f(x)|}\big]_{x^k}y^k$ for an arbitrary
non-constant function $f(x)$ which is not a polynomial of degree
two, then the spray {\bf G} in (\ref{ygj1}) is not
Finsler-metrizable. If {\bf G} in Theorem \ref{th2} has zero
Riemann curvature ($Ric=0$), then {\bf G} can be locally induced
by a Minkowski metric (a trivial case, Remark \ref{rem57} below).

For a spray of scalar curvature $R^i_{\ k}=R\delta^i_k-\tau_ky^i$
on a manifold $M$, the two quantities $R$ and $\tau_k$ are closely
related (see Proposition \ref{prop34} below). Now we consider the
following condition for $R$ and $\tau_k$,
 \be\label{ygj4}
 R_{.i}-2\tau_i=\omega_{i0},\ \ \ (\omega_{i0}:=\omega_{ir}y^r),
 \ee
where $\omega=\omega_{ij}dx^i\wedge dx^j$ is a 2-form on the
manifold $M$. A spray satisfying (\ref{ygj4}) is considered to be
some kind of weakly isotropic curvature. For a spray of scalar
curvature, the condition (\ref{ygj4}) is a special case of
$\chi_i=\omega_{i0}$ (see \cite{LS0}), where $\chi_i$ is called
the $\chi$-curvature originally defined in \cite{Shen1}. For a
spray satisfying (\ref{ygj4}), we have the following theorem.

 \begin{thm}\label{th3}
Let {\bf G} be an $n$-dimensional spray  of scalar curvature
$R^i_{\ k}=R\delta^i_k-\tau_ky^i$.
 \ben
  \item[{\rm (i)}] If  {\bf G} is induced by a Finsler
  metric $L$
 and satisfies (\ref{ygj4}). Then
  \ben
  \item[{\rm (ia)}] (\cite{FY1}) ($n\ge 3$) $L$ is of constant flag
curvature with $\omega=0$.
  \item[{\rm (ib)}] ($n=2$) For $L=F^2$, the flag curvature $\lambda$ of $L$ satisfies
 $\lambda''(\theta)+\epsilon I(\theta)\lambda'(\theta)=0$ on each tangent space, where
 $\epsilon=\pm 1$ is the sign of $L$, and $\theta$ is the Landsberg angle. Further, if $L$
 is a Riemann metric, or
 regular Finsler metric, then $L$ is of
isotropic flag curvature ($\lambda=\lambda(x)$) with $\omega=0$.
  \een
 \item[{\rm (ii)}] ($n=2$) If {\bf G} is induced by a Finsler
 metric $L>0$
 with constant main scalar, then {\bf G} satisfies (\ref{ygj4}) with $\omega$  not necessarily zero.
 \een
\end{thm}

Theorem \ref{th3}(i)   has essentially been proved in \cite{FY1}.
Starting from (\ref{ygj4}), we are also going to give a little
different version of proof from that in \cite{FY1}. If a spray
{\bf G} satisfying (\ref{ygj4}) with $\omega\ne 0$ and $n\ge 3$,
then the spray {\bf G} is not Finsler-metrizabl by Theorem
\ref{th3}(ia). If a two-dimensionalspray  ${\bf G}$ satisfying
(\ref{ygj4}) with $\omega\ne 0$, then ${\bf G}$ cannot be induced
by a Riemann metric or a regular Finsler metric by Theorem
\ref{th1}(ib), but Theorem \ref{th3}(ii) shows that possibly there
are such sprays which can be induced by a singular Finsler metric.

The problem still remains open to look for complete conditions for
 a spray of scalar curvature to be Finsler-metrizable. Only some
special conditions are concerned in Theorem \ref{th32} below.

\section{Preliminaries}\label{pre}

Let $M$ be an $n$-dimensional manifold. A conical region
$\mathcal{C}=\mathcal{C}(M)$ of $TM\setminus \{0\}$ means
$\mathcal{C}_x:=\mathcal{C}\cap T_xM\setminus \{0\}$ is conical
region for $x\in M$ ($\lambda y\in \mathcal{C}_x$ if
$\lambda>0,y\in \mathcal{C}_x$). A {\it spray} on $M$ is a smooth
vector field ${\bf G}$ on a conical region $\mathcal{C}$ of
 $TM\setminus \{0\}$ (an important case is $\mathcal{C}=TM\setminus \{0\}$) expressed in
a local coordinate system $(x^i,y^i)$ in $TM$ as follows
 $${\bf G}=y^i\frac{\pa}{\pa x^i}-2G^i\frac{\pa}{\pa y^i},$$
 where $G^i=G^i(x,y)$ are local functions satisfying
 $G^i(x,\lambda y)=\lambda^2G^i(x,y)$ for any constant $\lambda>0$.
The integral curves of {\bf G} projected onto $M$ are the
geodesics of {\bf G}.

The  Riemann curvature tensor $R^i_{\ k}$ of a given spray $G^i$
is defined by
 \be\label{y004}
  R^i_{\ k}:=2\pa_k G^i-y^j(\pa_j\dot{\pa_k}G^i)+2G^j(\dot{\pa_j}\dot{\pa_k}G^i)-(\dot{\pa_j}G^i)(\dot{\pa_k}G^j),
 \ee
where we define $\pa_k:=\pa/\pa x^k,\dot{\pa}_k:=\pa/\pa y^k$. The
trace of $R^i_{\ k}$ is called the Ricci curvature, $ Ric:=R^i_{\
i}$. A spray ${\bf G}$ is said to be {\it R-flat} if $R^i_{\
k}=0$. In \cite{Shen1}, Z. Shen defines a non-Riemannian quantity
called
 $\chi$-curvature $\chi=\chi_idx^i$ expressed as follows
  \be\label{j4}
 \chi_i:=2R^m_{\ i.m}+R^m_{\ m.i}.
  \ee
  Plugging (\ref{yg1}) into (\ref{j4}) yields
   \be\label{j5}
  \chi_i=(n+1)(R_{.i}-2\tau_i).
   \ee
   So a spray of scalar curvature is of isotropic curvature iff. it has vanishing
   $\chi$-curvature.

 A spray {\bf G} is called a Berwald spray if its Berwald curvature vanishes
 $G^i_{hjk}:=\dot{\pa}_h\dot{\pa}_j\dot{\pa}_kG^i=0$.
 A spray {\bf G} is said to be locally projectively flat if
 locally $G^i$ can be expressed as
$G^i=Py^i$, where $P$ is a positively homogeneous local function
of degree one.

In the calculation of some geometric quantities of a spray, it is
very convenient to use Berwald connection as a tool. For a spray
manifold $({\bf G}, M)$, Berwald connection is usually defined as
a linear connection on the pull-back $\pi^*\mathcal{C}$ ($\pi:
TM\rightarrow M$ the natural projection) over the base manifold
$\mathcal{C}$. The Berwald connection is defined by
  $$
  D(\pa_i)=(G^k_{ir}dx^r)\pa_k,\ \ \ \ \ (G^k_{ir}:=\dot{\pa}_r\dot{\pa}_iG^k),
  $$
 For a spray tensor $T=T_idx^i$ as an example,   the horizontal and vertical derivatives of $T$ with respect to Berwald
 connection are given by
  $$
 T_{i;j}=\delta_jT_i-T_rG^r_{ij},\ \ \ \ \ \ \
 T_{i.j}=\dot{\pa}_jT_i,\ \ \ \ (\delta_i:=\pa_i-G^r_i\dot{\pa}_r).
  $$
 The
 $hh$-curvature tensor $H^{\ i}_{j\ kl}$ of Berwald connection is
 defined by
 $$H^{\ i}_{j\ kl}:=\frac{1}{3}\big\{ R^i_{\ l.j.k}-(k/l)\big\},\ \
 H_{ij}:=H^{\ m}_{i\ jm},\ \ H_i:=\frac{1}{n-1}(nH_{0i}+H_{i0}),$$
where $T_{ij}-(i/j)$ means $T_{ij}-T_{ji}$,  and $T_0$ is defined
by $T_0:=T_ry^r$, as an
 example. For the Ricci identities and Bianchi identities of Berwald connection, readers
 can refer to \cite{AIM}.

 In this paper, we  define a Finsler metric $L(\ne 0)$ on a manifold $M$ as follows (cf. \cite{Shen6}):
 (i) for any $x\in M$, $L_x$ is defined on a conical region of  $T_xM\setminus \{0\}$
  and $L$ is $C^{\infty}$; (ii) $L$
 is positively homogeneous of degree two; (iii) the fundamental
metric tensor $g_{ij}:=(\frac{1}{2}L)_{y^iy^j}$ is non-degenerate.
 A Finsler metric
$L$ is said to be regular if additionally $L$ is defined on the
whole  $TM\setminus \{0\}$ and $(g_{ij})$ is positively definite.
Otherwise,
 $L$ is called singular. In general case, we
don't require that $L$ be regular. If a Finsler metric $L>0$, we
put $L=F^2$, and in this case, $F$ is also called a Finsler metric
being positively homogeneous of degree one.

 Any Finsler metric $L$ induces
a natural spray whose coefficients $G^i$ are given by
 $$
 G^i:=\frac{1}{4}g^{il}\big \{L_{x^ky^l}y^k-L_{x^l}\big
 \},
$$
where $(g^{ij})$ is the inverse of $(g_{ij})$. $L$ is said to be
of {\it scalar flag curvature} $K=K(x,y)$ if
 $$R^i_{\ k}=K(L\delta^i_k-y^iy_k),$$
 where $y_k:=(L/2)_{.k}=g_{km}y^m$. If $K_{.i}=0$, then $L$ is said to be of
 {\it isotropic flag curvature}, and in this case, $K$ is a
 constant if the dimension $n\ge 3$.

  A spray {\bf G} is (globally)
 Finsler-metrizable on $M$ (or on $\mathcal{C}(M)$) if there is a Finlser metric $L$ defined on
 a conical region $\mathcal{C}(M)$ and $L$ induces {\bf G}. A spray {\bf G} is
 locally Finsler-metrizable on $M$ if for each $x\in M$, there is
 a neighborhood $U$ of $x$ such that {\bf G} is
 Finsler-metrizable on $U$.

 Let $(M,F)$ be a
two-dimensional Finsler space with the Finsler metric $F$. We use
$\epsilon(=\pm1)$ to denote the sign of the determinant of the
metric matrix. We have
 \beq
&&g_{ij}=\ell_i\ell_j+\epsilon m_im_j=\ell_i\ell_j+h_{ij},\nonumber\\
&&\ell^i=y^i/F,\ \ \ \ (m^1,m^2)=\big(\sqrt{\epsilon
g}\big)^{-1}(-\ell_2,\ell_1),\ \ \  (g:=det(g_{ij})),\nonumber\\
&&Fm_{i.j}=-(\ell_i-\epsilon Im_i)m_j,\ \ \ \
FC_{ijk}=Im_im_jm_k,\label{cw9}
 \eeq
 where $\ell=(\ell^1,\ell^2),m=(m^1,m^2)$ is called the Berwald
 frame, $C_{ijk}$ is the Cartan tensor and $I$ is the main scalar.
 The system $F_{.i}=\ell_i,F\theta_{.i}=m_i$ is integrable. It
 defines the so-called  Landsberg angle $\theta$, which
is the arc-length parameter of the indicatrix $S_xM:=\{y\in
T_xM|F(x,y)=1\}$ with respect to the Riemann metric
$ds^2=g_{ij}dy^idy^j$ on the Minkowski plane $(M_x,F_x)$.

\begin{lem}\label{lem21}
 For a  positively homogeneous function
 $\lambda=\lambda(x,y)$ of degree zero, it satisfies on each
 tangent space,
 $$F\lambda_{.i}=\lambda'(\theta)m_i.$$
\end{lem}

\section{Spays of scalar curvature}

In this section, we will introduce some basic properties of sprays
 with scalar curvature, and the metrizability of such sprays
under certain  conditions.

\subsection{Some  basic formulas}

 For a spray of scalar curvature,
 $R$ and $\tau_k$ in  (\ref{yg1}) are related in the following formula (\ref{y8}).

 \begin{prop}\label{prop34}
 Let {\bf G} be an $n(\ge 3)$-dimensional spray  of scalar  curvature  $R^i_{\ k}=R\delta^i_k-\tau_ky^i$.
 Then there holds
  \be\label{y8}
   R_{.i;0}-3R_{;i}+\tau_{i;0}=0.
  \ee
  In particular, if {\bf G} is of isotropic curvature, then (\ref{y8}) becomes
  $\tau_{i;0}=R_{;i}$, or $\tau_{i;0}=\tau_{0;i}$.
\end{prop}

{\it Proof :} By a Bianchy identity of Berwald connection
  $$
 R^{i}_{\ jk;l}+R^{i}_{\ kl;j}+R^{i}_{\ lj;k}=0,
  $$
we have
 \be\label{y9}
 R^{m}_{\ k;m}+R^{m}_{\ km;0}-R^{m}_{\ m;k}=0.
 \ee
Since $R^i_{\ k}=R\delta^i_k-\tau_ky^i$, a direct computation
gives
 \be\label{y10}
R^{m}_{\ k;m}=R_{;k}-\tau_{k;0},\ \ \ R^{m}_{\
km;0}=\frac{1}{3}\big[(n-1)(R_{.k;0}+\tau_{k;0})+\tau_{k;0}-\tau_{m.k;0}y^m\big].
 \ee
By $\tau_0=R$, we obtain
 \be\label{y11}
 \tau_{m.k;0}y^m=R_{.k;0}-\tau_{k;0}.
 \ee
Now plugging (\ref{y10}) and (\ref{y11}) into (\ref{y9}) we have
 $$
 (n-2)(R_{.i;0}-3R_{;i}+\tau_{i;0})=0,
 $$
which gives the proof.   \qed

\

The isotropic case of Proposition \ref{prop34} is given by Z. Shen
in \cite{Shen2}. In Proposition \ref{prop34}, if {\bf G} is
Finsler-metrizable induced by a Finsler metric $L$ with the flag
curvature $\lambda$, we have $R=\lambda L, \tau_i=\lambda y_i$ and
then putting them into (\ref{y8}) gives a known formula
 \be\label{y12}
 L\lambda_{.i;0}+3\lambda_{;0}y_i-3L\lambda_{;i}=0.
 \ee
 If $L=F^2$ is of weakly isotropic flag curvature
$\lambda=3\theta/F+\sigma$, then (\ref{y12}) becomes
 $$
 \theta_{;i}-\theta_{.i;0}+(F\sigma_{;r}+2\theta_{;r})h^r_i=0,\ \
 \  (h^i_j:=\delta^i_j-\ell^i\ell_j).
 $$
For a two-dimensional spray {\bf G} in Proposition \ref{prop34},
(\ref{y8}) generally does not hold.  For example, if {\bf G} is a
two-dimensional spray induced by a Riemann metric $L$ of isotropic
Gauss curvature $\lambda=\lambda(x)$, then (\ref{y12}) reduces to
 $\lambda_{;i}L=\lambda_{;0}y_i$, which is impossible if
 $\lambda$ is not constant.

In Proposition \ref{prop34}, if {\bf G} is induced by a Finsler
metric $L$ of isotropic flag curvature $\lambda$, then we have
$\lambda=constant$ by (\ref{y12}), which is just the Schur's
Theorem. But for a general spray, we cannot conclude from
 Proposition \ref{prop34} that a spray of isotropic curvature in dimension $n\ge 3$ must be  of constant
curvature (see Examples \ref{ex73} and \ref{ex74} below).

 The following proposition gives a useful formula on a spray
 manifold of scalar curvature.

\begin{prop}\label{prop32}
 Let {\bf G} be a spray of scalar curvature $R^i_{\
k}=R\delta^i_k-\tau_ky^i$, and $T$ be a homogeneous scalar
function of degree $p$ satisfying $T_{;i}=0$. Then we have
 \be\label{yg61}
 RT_{.k}=pT\tau_k.
 \ee
 \ben
 \item[{\rm (i)}] If {\bf G} is of isotropic curvature with $R\ne
 0$, then there holds
  \be\label{yg62}
 T=c|R|^{\frac{1}{2}p}, \ \ \  (c=c(x)).
  \ee
\item[{\rm (ii)}] If {\bf G} is induced by a Finsler metric $L$ of
non-zero flag curvature, then there holds
 \be\label{yg63}
 T=c|L|^{\frac{1}{2}p},\ \ \ (c=constant).
 \ee
 \een
 \end{prop}

{\it Proof :}  By a Ricci identity and $R^i_{\
k}=R\delta^i_k-\tau_ky^i$, we have
 $$
 0=y^j(T_{;i;j}-T_{;j;i})=T_{.r}R^r_{\ i}
 =RT_{.k}-y^rT_{.r}\tau_k.
 $$
Then we obtain (\ref{yg61}) since $T$ is a homogeneous scalar
function of degree $p$.

 (i) Since {\bf G} is of isotropic curvature, we have
 $\tau_k=\frac{1}{2}R_{.i}$. Putting it into (\ref{yg61}) gives
 (\ref{yg62}) since
 $$
 \frac{T_{.i}}{T}=\frac{p}{2}\frac{R_{.i}}{R}\ \ \Big(\Leftrightarrow
 \big(\frac{T}{|R|^{\frac{1}{2}p}}\big)_{.i}=0\Big).
 $$

(ii) The original version is given in \cite{AIM}. Let
$\lambda\ne0$ be the flag curvature of $L$. Similarly as the proof
in (i), plugging $R=\lambda L$ and $\tau_i=\frac{1}{2}L_{.i}$ into
(\ref{yg61}) gives
$$T=c|L|^{\frac{1}{2}p},\ \ \  c=c(x).$$
 Then by $T_{;i}=0$ and  $L_{;i}=0$
we obtain $c_{;i}=0$, which means $c=constant$.  \qed

 \subsection{Some basic properties}
Let ${\bf G}$ be a spray of scalar  curvature satisfying
(\ref{yg1}). Then we have
 \be\label{y3}
H^{\ i}_{h\
jk}=\frac{1}{3}\big[R_{.j.h}\delta^i_k-\tau_{k.j.h}y^i-\tau_{k.j}\delta^i_h-\tau_{k.h}\delta^i_j
-(j/k)\big].
 \ee
  By (\ref{y3}) we obtain
 \be\label{y4}
 H_{ij}-H_{ji}=\frac{1}{3}(n+1)(\tau_{j.i}-\tau_{i.j}),
 \ee
 \be\label{y04}
 H_{0i}=\frac{1}{3}\big[(n-2)R_{.i}+(n+1)\tau_i\big],\ \ \ \
 H_{0i}=\frac{1}{3}\big[(2n-1)R_{.i}-(n+1)\tau_i\big]
 \ee
 \be\label{y5}
 H_{0i}-H_{i0}=-\frac{1}{3}(n+1)(R_{.i}-2\tau_i),
 \ee
 \be\label{y6}
 H_i=\frac{1}{3}(n+1)(R_{.i}+\tau_i).
 \ee

\begin{prop}\label{prop31}
 Let {\bf G} be a spray of scalar curvature
 satisfying (\ref{yg1}).
  \ben
   \item[{\rm (i)}]  If {\bf G} is of isotropic
 curvature, then
   $H_{ij}=H_{ji}$,  $H_{i0}=H_{0i}$ and $H_i$ is proportional to
   $\tau_i$.
   \item[{\rm (ii)}] If $H_{ij}=H_{ji}$, or $R_{i0}=H_{0i}$, or $H_i$ is proportional to
   $\tau_i$ with $R\ne0$, then {\bf G} is of isotropic
 curvature.
  \een

\end{prop}

{\it Proof :}  Firstly, we prove $R_{.i}=2\tau_i$ is equivalent to
$H_{ij}=H_{ji}$. If $R_{.i}=2\tau_i$, then it is easy to see that
$\tau_{i.j}=\tau_{j.i}$ and then by (\ref{y4}), we have
$H_{ij}=H_{ji}$. Conversely, if $H_{ij}=H_{ji}$, then by
(\ref{y4}) we have $\tau_{i.j}=\tau_{j.i}$. Differentiating
$R=\tau_0$ by $y^i$ yields
 $$R_{.i}=\tau_i+\tau_{m.i}y^m=\tau_i+\tau_{i.m}y^m=2\tau_i.$$

Secondly, it is clear from (\ref{y5}) that
 $R_{.i}=2\tau_i$ is equivalent to $H_{i0}=H_{0i}$.

Finally, if $R_{.i}=2\tau_i$, it is an obvious result that $H_i$
is proportional to $\tau_i$ by (\ref{y6}). Conversely, if $H_i$ is
proportional to $\tau_i$  and $R\ne 0$, then by (\ref{y6}) we get
$R_{.i}=\lambda \tau_i$ for some scalar function
$\lambda=\lambda(x,y)$. Contracting this by $y^i$ gives
$2R=\lambda \tau_0=\lambda R$. Since $R\ne0$ by assumption, we
have $\lambda=2$ and thus $R_{.i}=2\tau_i$.  \qed

\begin{prop}\label{prop33}
 A spray {\bf G} of scalar  curvature is R-flat if and only if
 $H_{ij}=0$, or $H_{i0}=0$, or $H_{0i}=0$ or $H_i=0$.
\end{prop}

{\it Proof :} The Riemann curvature tensor of {\bf G} satisfies
(\ref{yg1}). If $H_{ij}=0$, then $H_i=0$. So by (\ref{y6}) we have
$R_{.i}+\tau_i=0$. Contracting this by $y^i$ gives $3R=0$ and so
$R=0, \tau_i=0$. This shows that {\bf G} is R-flat. If $H_{0i}=0$
or $H_{i0}=0$, we have $(n-2)R_{.i}+(n+1)\tau_i=0$ or
$(2n-1)R_{.i}-(n+1)\tau_i=0$ by (\ref{y04}). Contracting either
one by $y^i$ gives $R=0$. Similarly we see that {\bf G} is R-flat.
\qed

\subsection{Metrizability under certain conditions}

For a spray of scalar curvature in general case, we have the
following result and Theorem \ref{th3} on metrizability.

\begin{thm}\label{th32}
Let {\bf G} be a spray of scalar curvature $R^i_{\
k}=R\delta^i_k-\tau_ky^i$.
 \ben
 \item[{\rm (i)}] If $R=\tau_k=0$, then {\bf G} is locally Finsler-metrizable.

 \item[{\rm (ii)}] If $R=0,\tau_k\ne 0$, then
 {\bf G} is not Finsler-metrizable.
 \item[{\rm (iii)}] If $R\ne 0$ and
  \be\label{yg15}
  \big(\frac{\tau_i}{R}\big)_{;j}\ne 0\ \ or\ \ \big(\frac{\tau_i}{R}\big)_{.j}\ne \big(\frac{\tau_j}{R}\big)_{.i},
  \ee
  then
 {\bf G} is not Finsler-metrizable.
 \een
\end{thm}

{\it Proof :} (i) is just Theorem \ref{th1}(i). For (ii), if {\bf
G} is induced by a Finsler metric $L$, then  $L$ is of scalar flag
curvature with $R=\lambda L, \tau_k=\lambda y_k$ for some scalar
function $\lambda=\lambda(x,y)$. So  if $R=0$, then $\tau_k=0$.
Thus {\bf G} is not Finsler-metrizable.

For (iii), if {\bf G} is induced by a Finsler metric $L$, then
similarly we have $R=\lambda L, \tau_k=\lambda y_k$. Thus we have
 $$
 \big(\frac{\tau_i}{R}\big)_{;j}=\big(\frac{y_i}{L}\big)_{;j}=0,\
 \ \
 \big(\frac{\tau_i}{R}\big)_{.j}=\big(\frac{y_i}{L}\big)_{.j}=\big(\frac{y_j}{L}\big)_{.i}=\big(\frac{\tau_j}{R}\big)_{.i}.
 $$
 So {\bf G} is not Finsler-metrizable.    \qed

\

In \cite{Shen6}, there are some two-dimensional sprays (of scalar
curvature) satisfying the condition of Theorem \ref{th32}(ii)
($R=0,\tau_k\ne 0$) (also see Example \ref{ex76} below).

\begin{rem}\label{rem36}
  Denote by $\mathcal{D}_{hol}( G)$ the holonomy algebra of
 a spray {\bf G} generated by the horizontal vector fields and their
 successive Lie brackets.  Define
 $$
 e_i:=(\delta^r_i-\frac{\tau_i}{R}y^r)\dot{\pa}_r, \ \
 (i=1,\cdots, n),
 $$
 where  {\bf G} has the dimension $n$ and is of scalar curvature $R^i_{\
k}=R\delta^i_k-\tau_ky^i$ with $R\ne0$.
 Then $\mathcal{V}:=span\{e_1,\cdots,e_n\}$ has the dimension $n-1$ and
 $\mathcal{V}\subset \mathcal{D}_{hol}( G)$ since  $Re_i=y^j[\delta_j,\delta_i]\in\mathcal{D}_{hol}(
 G)$. So $\mathcal{D}_{hol}( G)$ has the dimension at least $2n-1$. Further, we have
  $$
 [\delta_i,e_j]=
 G^r_{ij}e_r-\big(\frac{\tau_j}{R}\big)_{;i}Y,\ \ \ \ \
 [e_i,e_j]=\frac{\tau_i}{R}e_j-\frac{\tau_j}{R}e_i+\big[\big(\frac{\tau_i}{R}\big)_{.j}-
 \big(\frac{\tau_j}{R}\big)_{.i}\big]Y,\ \ \ \ (Y:=y^i\dot{\pa}_i).
  $$
 Note that $Y$ is not in $\mathcal{V}$. Here we can give a
 different proof of
   Theorem \ref{th32}(iii): if (\ref{yg15}) holds, then we
   have $Y\in \mathcal{D}_{hol}(G)$ and thus $0=YL=2L$ (a
   contradiction). On the other hand, even if
   \be\label{cw16}
 \big(\frac{\tau_i}{R}\big)_{;j}= 0\ \ and\ \ \big(\frac{\tau_i}{R}\big)_{.j}= \big(\frac{\tau_j}{R}\big)_{.i},
   \ee
   we still do not make sure whether $Y$ is in
    $\mathcal{D}_{hol}( G)$ or not.

    It is easy to check that a spray  of cosntant  curvature
 ($R\ne 0$) always satisfies (\ref{cw16}) (see (\ref{ycw65}) below). The condition (\ref{cw16}) is a necessary but not sufficient for
 {\bf G} to be Finsler-metrizable in Theorem \ref{th32}(iii). See Example \ref{ex74} below based on Theorem \ref{th2}.
\end{rem}

{\it Proof of Theorem \ref{th3}  :} Let {\bf G} be induced by the
Finsler metric $L$. Here we provide a version of proof of Theorem
\ref{th3}(ia) for $L>0$. Put $R=\lambda L, \tau_i=\lambda
y_i=\lambda F\ell_i$, where $F$ is given by $F^2=L$ and $\lambda$
is the flag curvature. By (\ref{ygj4}) we have
$(R_{.i}-2\tau_i)_{.k}+(i/k)=0$, which just is
 \be\label{y18}
 F\lambda_{.i.k}+\lambda_{.k}\ell_i+\lambda_{.i}\ell_k=0.
  \ee
 Differentiating (\ref{y18}) by
 $y^j$ gives
 \be\label{y17}
F\lambda_{.i.j.k}+\lambda_{.i.k}\ell_j+\lambda_{.j.k}\ell_i+F^{-1}\lambda_{.k}h_{ij}+
\lambda_{.i.j}\ell_k+F^{-1}\lambda_{.i}h_{jk}=0.
 \ee
Interchanging $i,j$ in (\ref{y17}) and making a subtraction, we
obtain
 $
 h_{jk}\lambda_{.i}-h_{ik}\lambda_{.j}=0.
 $
Contracting this by $g^{jk}$ gives $(n-2)\lambda_{.i}=0$.

(i)(ia) If $n>2$, then we have $\lambda_{.i}=0$. So $\lambda$ is
constant by Schur's Theorem.

(i)(ib) If $n=2$, let $(\ell,m)$ be the Berwald frame and $\theta$
be the Landsberg angle.  By Lemma \ref{lem21}, we put
 \be\label{g16}
 F\lambda_{.i}=\eta m_i, \ \ \ \
 \eta(\theta)=\lambda'(\theta).
 \ee
 Differentiating (\ref{g16}) by $y^j$ and using (\ref{cw9}) we obtain
  \beq
 \ell_j\lambda_{.i}+F\lambda_{.i.j}\hspace{-0.6cm}&&=\eta_{.j}m_i-\eta
 F^{-1}(\ell_i-\epsilon Im_i)m_j,\nonumber\\
\ell_i\lambda_{.j}+\ell_j\lambda_{.i}+F\lambda_{.i.j}
\hspace{-0.6cm}&&=\eta_{.j}m_i+\epsilon\eta
F^{-1}Im_im_j.\label{g017}
  \eeq
  Then by (\ref{y18}) we have
  $$
 F\eta_{.j}+\epsilon\eta Im_j=0,\ \ \ {\rm or}\ \ \  \eta'+\epsilon\eta
 I=0,\ \ \ {\rm or}\ \ \  \lambda''+\epsilon\lambda' I=0.
  $$
 If $F$ is
 regular ($\epsilon=1$), then $T_xM$ is compact. There is a $\theta_0$ such that
$\lambda'(\theta_0)=0$. So we get
$\eta(\theta)=\lambda'(\theta)\equiv 0$. Thus we have
$\lambda_{.i}=0$ from (\ref{g16}).

(ii) It is  known that all  Finsler metrics with isotropic main
scalar $I=I(x)$ on a two-dimensional manifold are divided into the
following three classes (\cite{AIM}):
 \beq
   &&L=\beta^{2s}\gamma^{2(1-s)},\quad (s=s(x)\neq 0,\quad s(x)\neq 1),\label{g26}\\
  &&L=\beta^2e^{\frac{2\gamma}{\beta}},\label{g27}\\
  &&L=(\beta^2+\gamma^2)e^{2r\cdot arctan(\frac{\beta}{\gamma})},\quad r=r(x), \label{g28}
\eeq
 where $\beta=p_i(x)y^i$ and $\gamma=q_i(x)y^i$ are two
independent 1-forms. The main scalar $I=I(x)$ is given
respectively by
 $$
  \epsilon I^2=\frac{(2s-1)^2}{s(s-1)},\ \ \
  I^2=4,\ \ \ \
  I^2=\frac{4r^2}{1+r^2}.
$$
In the following we assume that the Finsler metric $L$ is of
constant main scalar, that is, $s$ in (\ref{g26}) and $r$ in
(\ref{g28}) are constant. For the convenience of computation, we
may put $\beta=py^1$ and $\gamma=qy^2$ under certain local
coordinate system, where $p=p(x^1,x^2)$ and $q=q(x^1,x^2)$ are
scalar functions. By a direct computation, we can obtain the flag
curvature $\lambda$ of the Finsler metric given by
(\ref{g26})-(\ref{g28}) respectively (see \cite{YC}). Now we can
obtain (\ref{ygj4}) from $R=\lambda L$ and $\tau_k=\lambda y_k$.
If $L$ is given by (\ref{g26}), then (\ref{ygj4}) holds with
 $$
 \omega_{12}=\frac{1-2s}{s(1-s)}\frac{(p^2qq_{12}-pq^2p_{12}+q^2p_1p_2-p^2q_1q_2)s-p^2qq_{12}+p^2q_1q_2}{p^2q^2}.
 $$
If $L$ is given by (\ref{g27}), then (\ref{ygj4}) holds with
 $$
 \omega_{12}=-\frac{2(p^2qp_{22}+pq^2p_{12}-p^2qq_{12}+p^2q_1q_2-q^2p_1p_2-p^2p_2q_2)}{p^2q^2}.
 $$
If $L$ is given by (\ref{g28}), then (\ref{ygj4}) holds with
 $$
 \omega_{12}=\frac{2r}{1+r^2}\frac{(p^2qq_{12}-pq^2p_{12}-p^2q_1q_2+q^2p_1p_2)r+p^2qp_{22}+pq^2q_{11}-q^2p_1q_1-p^2p_2q_2}{p^2q^2}.
 $$
So generally, $\omega_{12}$ in the above is not zero.    \qed

\section{Sprays of isotropic curvature}
 In this section, we are going to prove Theorem \ref{th1} and give
 a metrizibility result as an application of Theorem \ref{th1}. For this, we first introduce a simple necessary and sufficient
 condition for a general spray to be metrizable (see Lemma \ref{lem043}).

 In the
following lemmas, the horizontal and vertical covariant derivative
is taken with respect to the given spray {\bf G}.

\begin{lem}\label{lem041}
Let {\bf G} be a spray and $L$ be a Finsler function. Then  we
have
 $$
L_{.i;0}-L_{;i}=0\Longleftrightarrow G^i=\frac{1}{4}g^{il}\big
\{L_{x^ky^l}y^k-L_{x^l}\big
 \}.
 $$
\end{lem}

{\it Proof :} It follows  from
 $
L_{.i;0}-L_{;i}=-2G^rL_{.r.i}+L_{x^ry^i}y^r-L_{x^i}.
 $
\qed

\begin{lem}\label{lem042}
Let {\bf G} be a spray and $L$ be a Finsler function. Then we have
 $$
L_{.0}=2L,\ \ \ L_{;i}=0\ \Longleftrightarrow \  L_{.0}=2L,\ \ \
L_{.i;0}-L_{;i}=0.
 $$
\end{lem}

{\it Proof :} We only need to prove "$\Longleftarrow$". By
$L_{.i;0}=L_{;i}$ we have
 $L_{;0.i}-L_{;i}=L_{;i}$, or $L_{;0.i}=2L_{;i}$. Further,
 $L_{.i;0}=L_{;i}$ implies $L_{.0;0}=L_{;0}$. So by $L_{.0}=2L$ we
 obtain $L_{;0}=0$. Thus it follows from $L_{;0.i}=2L_{;i}$ that
 $L_{;i}=0$.   \qed

\begin{lem}\label{lem043}
A spray {\bf G} is Finsler-metrizable if and only if there is  a
Finsler function $L$ satisfies $L_{;i}=0$. In this case, {\bf G}
is  induced by $L$.
\end{lem}

  \subsection{Formal integrability}

   To prove Theorem \ref{th1}(i), we need a theory  on
   Spencer's technique of formal integrability for linear partial
   differential systems. Here we only give some basic notions for this theory
    and more details are refereed to \cite{GM,BCG}.

 Let $B$ be a vector bundle over an $n$-dimensional manifold $M$, and denote by
 $J_kB$  the bundles of $k$th order jets of the sections of $B$.
 For two vector bundles $B_1, B_2$ over $M$, consider $P:
 Sec(B_1)\rightarrow Sec(B_2)$, which is a linear partial
   differential operator of order $k$. $P$ can be identified with
   a map $p_0(P): J_kB_1\rightarrow B_2$, a morphism of vector
   bundles over $M$. We also denote by $p_l(P): J_{k+l}B_1\rightarrow
   J_lB_2$ the morphisms of vector bundles over $M$, which is
   called the $l$th order jet prolongation of $P$. Let
   $R_{k+l,x}(P):=Ker p_l(P)_x$ be the space of $(k+l)$th order
   formal solutions of $P$ at a point $x\in M$. The operator $P$
   is said to be formally integrable at $x\in M$, if $R_{k+l}(P)$
   is a vector bundle for all $l\ge 0$ and the projection
   $\pi_{k+l,x}(P):R_{k+l,x}(P)\rightarrow R_{k+l-1,x}(P)$ is onto
   for all $l\ge 1$.

 Let $\sigma_k(P):S^k(T^*M)\otimes B_1\rightarrow B_2$
 be the symbol of $P$, which is defined by the highest  order term
 of $P$, and let $\sigma_{k+l}(P):S^{k+l}(T^*M)\otimes B_1\rightarrow S^l(T^*M)\otimes
 B_2$ be the symbol of the $l$th order prolongation of $P$. Define
  \beqn
 g_{k,x}(P):&&\hspace{-0.6cm}=Ker\
 \sigma_{k,x}(P),\\
  g_{k,x}(P)_{e_1\ldots
 e_j}:&&\hspace{-0.6cm}=\{A\in g_{k,x}(P)|i_{e_1}A=\ldots =i_{e_j}A=0\},\ \ \ 1\le j\le
 n,
  \eeqn
where $\{e_1,\cdots, e_n\}$ is a basis of $T_xM$. Such a basis is
said to be quasi-regular if it satisfies
 $$
 dim\ g_{k+1,x}(P)= dim\ g_{k,x}(P)+\sum^n_{j=1}dim\
 g_{k,x}(P)_{e_1\cdots e_j}.
 $$
The symbol $\sigma_k(P)$ is said to be involutive at $x\in M$ if
there exists a quasi-regular  basis of $T_xM$. For the proof of
Theorem \ref{th1}(i), we need the following theorem and lemma
(\cite{GM}).

\begin{thm} (Cartan-Kahler) \label{th41} Let $P$ be a regular linear partial
differential operator of order $k$. If
$\pi_{k+1,x}(P):R_{k+1,x}(P)\rightarrow R_{k,x}(P)$ is onto and
the symbol $\sigma_k(P)$ is involutive, then $P$ is formally
integrable.
\end{thm}

\begin{lem}\label{lem42}
For two vector bundles $B_1, B_2$ over $M$, let $P:
Sec(B_1)\rightarrow Sec(B_2)$ be a regular linear partial
differential operator of order $k$. Then
$\pi_{k+1,x}(P):R_{k+1,x}(P)\rightarrow R_{k,x}(P)$ is onto iff.
 $$
 P(s)_x=0\Longrightarrow (DP(s))_x=\sigma_{k+1}(P)(A):
 $$
 for some $A\in S^{k+1}(T^*_xM)\otimes B_1$, where $D$ is an arbitrary linear
 connection of the bundle $B_2$ over the base manifold $M$.
\end{lem}

\subsection{Proof of Theorem \ref{th1}(i)}
 Let $T^*_vTM$ denote the subbundle of $T^*TM$, in which, if
$\omega\in T^*_vTM$, then $\omega$ can be written locally as
$\omega=\omega_i(x,y)dx^i$, and $\wedge^2T^*_vTM$ the subbundle of
$T^*TM$, every element $\omega$ of which is locally in the form
$\omega=\omega_{ij}(x,y)dx^i\wedge dx^j$ with
$\omega_{ij}=-\omega_{ji}$. Besides, we designate $S^k(T^*TM)$ as
the bundle of symmetric $k$-forms over $TM$. For an
$n$-dimensional manifold $M$, let $B_1,B_2$ be two vector bundles
over $TM$ with
 $$
 B_1:= T^*_vTM,\ \ \ B_2:=T^*_vTM\oplus\wedge^2
T^*_vTM\oplus(T^*_vTM\otimes T^*_vTM).
 $$
We define a linear partial differential operator
$P:Sec(B_1)\rightarrow Sec(B_2)$ in component form as follows
 \be\label{yg038}
 P(\theta_i)=(\theta_{i.0},\theta_{i.j}-\theta_{j.i},
 \theta_{i;j}).
 \ee

\begin{lem}\label{lem43}
 For the operator $P$ in (\ref{yg038}), the symbol $\sigma_1(P)$ is involutive.
\end{lem}

\noindent{\it Proof :} We are going to prove that
$\{\dot{\pa}_1,\cdots,\dot{\pa}_n,\delta_1,\cdots,\delta_n\}$ is a
quasi-regular basis of $P$.

 By definition, for
$A=(A_{ji},A_{\underline{j}i})\ (=A_{ji}dx^j\otimes
dx^i+A_{\underline{j}i}\delta y^j\otimes dx^i)\in T^*TM\otimes
B_1$, we have
 $$
\sigma_1(P)A=(A_{\underline{0}i},A_{\underline{j}i}-A_{\underline{i}j},A_{ji})\in
B_2.
 $$
Assume $\sigma_1(P)(A)=0$. Then for the computation of $dim
(g_1(P))$, we see that $A_{\underline{0}i}=0$ and
$A_{\underline{i}j}=A_{\underline{j}i}$ together contribute the
number $(n-1)n/2$, and $A_{ij}=0$ gives $0$. Therefore, we obtain
 \be\label{yg38}
 dim(g_1(P))=\frac{(n-1)n}{2}.
 \ee

 Now with respect to the basis $\{dx^i,\delta y^i\}$, an element
 $B\in S^2(T^*TM)\otimes B_1$ can be expressed as
  \beqn
 &&\hspace{0.8cm} B=(B_{ijk},\ B_{\underline{i}jk},\ B_{i\underline{j}k},\ B_{\underline{i}\underline{j}k})
\\
 &&\big(B_{ijk}=B_{jik},B_{\underline{i}jk}=B_{j\underline{i}k},B_{\underline{i}\underline{j}k}=B_{\underline{j}\underline{i}k}\big).
 \eeqn
By definition, we have
 $$
 \sigma_2(P)B=(B_{i\underline{0}k},B_{\underline{i}\underline{0}k};B_{i\underline{j}k}-B_{i\underline{k}j},
B_{\underline{i}\underline{j}k}-B_{\underline{i}\underline{k}j};B_{ijk},B_{\underline{i}jk}).
 $$
Assume $\sigma_2(P)(B)=0$. Then $B_{ijk}$, $B_{i\underline{j}k}$
and $B_{\underline{i}jk}$ gives 0 to $dim(g_2(P))$.  By the fact
that
 $B_{\underline{i}\underline{j}k}$ is symmetric in $i,j,k$ satisfying
additional condition $B_{\underline{i}\underline{0}k}=0$, we
obtain the number
 $$
\sum^{n-2}_{k=0}\frac{(n-k-1)(n-k)}{2}.
 $$
to $dim(g_2(P))$. So altogether, we have
 \be\label{yg39}
dim(g_2(P))=\sum^{n-2}_{k=0}\frac{(n-k-1)(n-k)}{2}.
 \ee

Next we verify  under the basis
$\{\dot{\pa}_1,\cdots,\dot{\pa}_n,\delta_1,\cdots,\delta_n\}$ at a
point $(x,y)\in TM$,
 \be\label{yg40}
 dim(g_2(P))=dim\ g_{1}(P)+\sum^n_{j=1}dim\
 g_{1}(P)_{\dot{\pa}_1\cdots \dot{\pa}_j}+\sum^n_{j=1}dim\
 g_{1}(P)_{\dot{\pa}_1\cdots \dot{\pa}_n\delta_1\cdots \delta_j}.
 \ee
By a direct computation, we see
 \beqn
 dim\ g_{1}(P)_{\dot{\pa}_1\cdots \dot{\pa}_j}=\frac{(n-j-1)(n-j)}{2},\
 \ \ \
 dim\ g_{1}(P)_{\dot{\pa}_1\cdots \dot{\pa}_n\delta_1\cdots \delta_j}
 =0,\ \ (1\le j\le
 n).
 \eeqn
Plugging them into (\ref{yg40}) and using (\ref{yg38}),
(\ref{yg39}), we see that both sides of (\ref{yg40}) are equal.
Therefore, $\sigma_1(P)$ is involutive.  \qed

\begin{lem}\label{lem44}
 For the operator $P$ in (\ref{yg038}), a first-order solution of
 $P(\theta)=0$ can be lifted into a second-order solution iff.
  \be\label{yg42}
 \theta_rH^{\ r}_{i\ jk}+\theta_{i.r}R^r_{\ jk}=(\theta_{r}R^r_{\ jk})_{.i}=0.
  \ee
\end{lem}

{\it Proof :} If $P(\theta)=0$, then we have
$\theta_{i;j}=\theta_{j;i}$. Then by a Ricci identity and
$\theta_{i.j}=\theta_{j.i}$, it gives (\ref{yg42}):
 $$
  0=\theta_{i;j;k}-\theta_{i;k;j}=-\theta_rH^{\
  r}_{i\
  jk}-\theta_{i.r}R^r_{jk}=-(\theta_{r}R^r_{\ jk})_{.i}.
  $$

Conversely, suppose that (\ref{yg42}) holds for a $\theta$
satisfying $P(\theta)=0$ with $\theta_0\ne 0$. Let $D$ be the
Berwald connection of the bundle $\pi^*TM$ over the base manifold
$TTM$. Then $D$ can be naturally extended to the bundle $B_2$. We
use Lemma \ref{lem42}. In component form, we have
 \be\label{yg43}
 DP(\theta_i)=(\theta_{i.0;j},\ \theta_{i.0.j},\ \
 \theta_{i.j;k}-\theta_{j.i;k},\ \theta_{i.j.k}-\theta_{j.i.k},\ \
 \theta_{i;j;k},\ \theta_{i;j.k}),
 \ee
where $\theta$ ($\theta_0\ne 0$) satisfies $P(\theta)=0$ at a
point $w=(x,y)\in TM$, that is, at the point $w$ there holds
 \be\label{yg44}
 \theta_{i.0}=0,\ \ \ \ \theta_{i.j}=\theta_{j.i},\ \ \ \
 \theta_{i;j}=0.
 \ee
Next we are going to prove that there is an $A_{\alpha\beta i}\in
S^2(T^*_wTM)\otimes B_1$ satisfying
 \be\label{yg45}
DP(\theta_i)=\sigma_2(A_{\alpha\beta
i})=(A_{\alpha\underline{r}i}y^r,\ \
A_{\alpha\underline{j}i}-A_{\alpha\underline{i}j},\ \ A_{\alpha
ji}).
 \ee
By (\ref{yg43}), we see that (\ref{yg45}) is equivalent to
 \beqn
 &&\text{\ding{172}}\ A_{j\underline{r}i}y^r=\theta_{i.r;j}y^r,\hspace{3.1cm}
  \text{\ding{173}}\
  A_{\underline{j}\underline{r}i}y^r=\theta_{i.0.j}=\theta_{i.r.j}y^r+\theta_{i.j},\\
  &&\text{\ding{174}}\
  A_{k\underline{j}i}-A_{k\underline{i}j}=\theta_{i.j;k}-\theta_{j.i;k},\hspace{1.5cm}
 \text{\ding{175}}\
 A_{\underline{k}\underline{j}i}-A_{\underline{k}\underline{i}j}=\theta_{i.j.k}-\theta_{j.i.k},\\
 &&\text{\ding{176}}\
 A_{kji}=\theta_{i;j;k},\hspace{3.75cm}
\text{\ding{177}}\ A_{\underline{k}ji}=\theta_{i;j.k}.
 \eeqn
 In the following, we will construct $A_{\alpha\beta i}$ which
 satisfies the above six relations. Put
 \beqn
&&A_{\underline{j}\underline{k}i}=\theta_{i.k.j}+\theta_0^{-1}\big[\theta_{i.j}\theta_k+(i,j,k)\big],\\
&&A_{j\underline{k}i}=\theta_{i.k;j}-\theta_rG^r_{ijk},\ \ \ \ \
A_{\underline{k}ji}=\theta_{i;j.k},\\
&&A_{kji}=\theta_{i;j;k}.
 \eeqn
Firstly, $A_{\underline{j}\underline{k}i}$ is symmetric in $j,k$
by $\theta_{i.j}=\theta_{j.i}$ in (\ref{yg44}) and it also
satisfies \ding{173}\ding{175} by $\theta_{i.j}=\theta_{j.i},
\theta_{i.0}=0$ in (\ref{yg44}). Next, by a Ricci identity of
Berwald connection, we see that
 $$
A_{\underline{k}ji}-A_{j\underline{k}i}=\theta_{i;j.k}-\theta_{i.k;j}+\theta_rG^r_{ijk}=0,
 $$
 which gives $A_{\underline{k}ji}=A_{j\underline{k}i}$. It is also
 clear that $A_{\underline{k}ji}$ and $A_{j\underline{k}i}$
 satisfy \ding{172}\ding{174}\ding{177}. Finally, $A_{kji}$ is
 symmetric in $j,k$, and satisfies \ding{176} by a Ricci identity of Berwald
 connection and (\ref{yg42}). So (\ref{yg45}) holds. This finishes
 the proof of the lemma.   \qed

\

Now in Theorem \ref{th1}(i), we have $R=0$. So (\ref{yg42})
automatically holds.  It follows from Lemmas \ref{lem43},
\ref{lem44} and then Theorem \ref{th41}, the operator $P$ is
formally integrable, that is, for each point $u_0:=(x^i_0,y^i_0)$,
there exist a neighborhood $U$ of $u_0$ and a analytic $\theta$
defined on $U$ such that $P(\theta)=0$.

Under the basis $(\delta_i,\dot{\pa}_i)$, the local coordinate of
$J_1B_1$ is expressed as
$(x^i,y^i,\theta_i,\theta_{i\underline{j}}, \theta_{ij})$. An
initial data $(x^i_0,y^i_0,\theta_i^0,\theta^0_{i\underline{j}},
\theta^0_{ij})$ satisfied by the operator $P$ means
$\theta_{i\underline{0}}^0=0,\theta^0_{i\underline{j}}=\theta^0_{j\underline{i}},\theta^0_{ij}-\theta_r^0G^r_i=0$.
Further, we let the initial data satisfy $\theta_i^0y^i_0>0$ and
$Rank(\theta^0_{i\underline{j}})=n-1$.

Now for the above analytic solution $\theta$ of $P$ which is
defined on a neighborhood $U$ of $u_0$ and satisfies the above
initial data, we obtain a local  metric $F:=\theta_0$ defined on
$U$. To prove that $F$ is a Finsler metric, we need the following
lemma which can be proved by an elementary discussion in linear
algebra.

\begin{lem}
 Let $F$ be positively homogeneous of degree one with $F(y)\ne 0$
 at a point $y$. Then $g_{ij}:=\frac{1}{2}(F^2)_{.i.j}$ is
 non-degenerate at $y$ iff. $Rank(F_{.i.j})=n-1$ at $y$.
\end{lem}

Now for $F=\theta_0$, we have $F_{.i.j}=\theta_{i.j}$. Then at
$u_0$, we have
$Rank(F_{.i.j})=Rank(\theta^0_{i\underline{j}})=n-1$. So by the
above lemma, $g_{ij}$ is non-degenerate at $u_0$. By continuity,
$g_{ij}$ is non-degenerate in $U$ (when it is small enough). Thus
we obtain a Finsler metric $F$ with each $F_x$ defined on the
conical region formed by $y=0$ and  $\{y|(x,y)\in U\}$. Further,
$F$ satisfies $F_{;i}=0$, which means that the spray {\bf G} in
Theorem \ref{th1} is induced by $F$ by Lemma \ref{lem043}. This
completes the proof of Theorem \ref{th1}(i). \qed

\begin{rem}
 The idea of the proof of Theorem \ref{th1}(i) can be referred
to
 that in \cite{BM} for the formal integrability of the operator
 $P_1(\theta):=(\theta_{i.0},\theta_{i.j}-\theta_{j.i},\delta_i\theta_j-\delta_j\theta_i)$.
 On the other hand, a suitable change of the proof of Theorem \ref{th1}(i)
 can give the proof for the formal integrability of the operator
 $P_1$ in \cite{BM} (where actually we can redefine $P_1$ as $\bar{P}_1:=
 (\theta_{i.0},\theta_{i.j}-\theta_{j.i},\theta_{j;i}-\theta_{i;j})$).
 Besides, we may also consider  the system $F_{.0}=F, \ F_{;i}=0$ for a
 possible proof (cf. \cite{Mu} for a more general discussion).
\end{rem}

\subsection{Proof of Theorem \ref{th1}(ii) and (iii)}

 \hspace{0.35cm}  (ii) Assume that $R\ne 0$ is not a Finsler metric.
If {\bf G} is Finsler-metrizable induced by a Finsler metric $F$,
then $F$ is of isotropic curvature $\lambda\ne 0$. By $R=\lambda
F^2$, we see that $R$ is a Finsler metric, which gives a
contradiction.

(iii) Assume that $R$ is a Finsler metric. If {\bf G} is
Finsler-metrizable induced by a Finsler metric $L$, then $L$ is of
isotropic curvature $\lambda=\lambda(x)\ne 0$. By $R=\lambda L$,
we have $R_{;i}=\lambda_{;i}L$. Thus we obtain
 $$
R_{;i}=\frac{\lambda_{;i}}{\lambda}R=R(\ln|\lambda|)_{;i}.
 $$
Let $\omega_i:=(\ln|\lambda|)_{;i}$. Then $\omega$ is closed and
$R_{;i}=R\omega_i$.  Conversely, if $R_{;i}=R\omega_i$ for some
closed 1-form $\omega=\omega_i(x)dx^i$, then locally there is a
scalar function $\lambda=\lambda(x)\ne 0$ such that
$\omega_i=(\ln|\lambda|)_{;i}$. It is easy to check that
$R_{;i}=R(\ln|\lambda|)_{;i}$ is equivalent to
$(R/\lambda)_{;i}=0$. Therefore, {\bf G} is Finsler-metrizable
induced by the Finsler metric $L:=R/\lambda$ by Lemma
\ref{lem043}. By $R=\lambda L$, we see that $L$ is of isotropic
flag curvature $\lambda$. If $\omega=0$, we
      may choose $\lambda=1$, and then the Finlser metric
      $L=R$ is of constant flag curvature $\lambda=1$. If $n\ge3$, then the Finsler metric $L$ is of constant flag
curvature $\lambda$ by Schur's theorem, which gives $R_{;i}=0$.
\qed

\subsection{A metrizability result}
As an application of Theorem \ref{th1}, we show the following
theorem.

 \begin{thm}\label{th410}
Let $G^i$ be the spray of a Finsler metric $F$ of constant flag
curvature $\lambda$ and $\bar{G}^i$ be a spray defined by
$\bar{G}^i=G^i+cFy^i$ for a constant $c$. Then ${\bf \bar{G}}$ is
(locally) Finsler-metrizable iff. $\lambda=-c^2$ or $c=0$. When
$\lambda=-c^2$, $\bar{G}^i$ is locally induced by a Finsler metric
of zero flag curvature.
\end{thm}

{\it Proof :} The Riemann curvature $R^i_{\ k}$ of {\bf G} is
given by $$R^i_{\ k}=\lambda(F^2\delta^i_k-FF_{.k}y^i).$$ Then by
a direct computation, the Riemann curvature $\bar{R}^i_{\ k}$ of
${\bf \bar{G}}$ is given by
 $$
\bar{R}^i_{\ k}=\bar{R}\delta^i_k-\bar{\tau}_ky^i,\ \ \ \
\big(\bar{R}:=(\lambda+c^2)F^2,\ \ \
\bar{\tau}_k:=(\lambda+c^2)FF_{.k}\big).
 $$
So ${\bf \bar{G}}$ is of isotropic curvature since
$\bar{R}_{.i}=2\bar{\tau}_i$. Further, we have
 $$
 \bar{R}_{\bar{;}i}=\bar{R}_{;i}+\bar{R}_{.r}(cF_{.i}y^r+cF\delta^r_i)=4c(\lambda+c^2)F^2F_{.i}.
 $$
Assume that ${\bf \bar{G}}$ is Finsler-metrizable. If
$c(\lambda+c^2)\ne 0$, then $\bar{R}$ is a Finsler metric. By
Theorem \ref{th1}(iii), we have
$\bar{R}_{\bar{;}i}=\bar{R}\omega_i$ for some closed 1-form
$\omega=\omega_i(x)dx^i$. But clearly this does not hold.
Therefore, we have $c=0$ or $\lambda=-c^2$. Conversely, if $c=0$,
then ${\bf \bar{G}}$  is induced by $F$. If $\lambda+c^2= 0$, then
${\bf \bar{G}}$ has zero Riemann curvature. So ${\bf \bar{G}}$ is
(locally) Finsler-metrizable by Theorem \ref{th1}(i).\qed

\

Theorem \ref{th410} is a  generalization of a result in
\cite{Yang1}, where we have an additonal condition that $F$ is
projectively flat.

\section{Sprays of constant curvature}

In this part, we introduce a new notion: a spray of {\it constant
curvature}, which is a generalization of a Finsler metric of
constant flag curvature. For this new notion, some basic
properties for Finsler metrics still remain unchanged  for sprays
(see Theorem \ref{th01} below).

 \begin{Def}
A spray {\bf G}  of scalar curvature $R^i_{\
k}=R\delta^i_k-\tau_ky^i$ is said to be of {\it constant
curvature} if  $\tau_{i;j}=0$.
 \end{Def}

The following theorem gives some basic properties for sprays of
constant curvature.

\begin{thm}\label{th01}
A  spray  has the following properties on constant curvature:
  \ben
 \item[{\rm (i)}] A spray of scalar curvature  is of constant curvature iff.
  its Riemann curvature is zero or its Ricci curvature $Ric$ satisfies $Ric_{;i}=0 (Ric\ne0)$.
 \item[{\rm (ii)}] A spray of  constant
 curvature must be of isotropic curvature.
 \item[{\rm (iii)}]
 A Finsler metric  is of constant flag
curvature iff. its spray is of  constant curvature. \item[{\rm
(iv)}]  An $n$-dimensional spray of isotropic curvature is not
necessarily of constant curvature
 even for $n\ge 3$.
 \een
 \end{thm}

A Finsler metric has the same conclusions as shown for sprays in
Theorem \ref{th01}(i)(ii).  Meanwhile,  Theorem \ref{th01}(iv)
shows a different property of sprays from that of Finsler metrics.

To prove Theorem \ref{th01}, we first show the following lemma.

\begin{lem}
Let {\bf G} be a spray of scalar curvature $R^i_{\
k}=R\delta^i_k-\tau_ky^i$.  Then we have
  \beq
  &&R_{;i}=0\ (R\ne0)\ \Longrightarrow \  R_{.i}=2\tau_i,\label{ycw64}\\
&&\tau_{i;k}=0\ \Longrightarrow \ R=\tau_k=0\ \ or \ \ R_{;i}=0\
(R\ne0),\label{ycw65}\\
&&R_{;i}=0\ (R\ne0)\ \Longrightarrow \  \tau_{i;k}=0.\label{ycw66}
  \eeq
 \end{lem}

{\it Proof :}  Assume  $R_{;i}=0 (R\ne 0)$.
 We have $RR_{.i}=2R\tau_i$ from (\ref{yg61}), where we have put $T=R$ with $p=2$ in (\ref{yg61}). This gives
 $R_{.i}=2\tau_i$ since $R\ne0$, which gives the proof of
 (\ref{ycw64}).

Assume $\tau_{i;k}=0$. We have $R_{;i}=0$ since $\tau_0=R$ and
then $R_{;i}=\tau_{0;i}=\tau_{m;i}y^m=0$. Now we prove that if $R=
0$, then $\tau_i=0$. By a
 Ricci identity we obtain
  \be\label{yr64}
 0=y^j(\tau_{i;j;k}-\tau_{i;k;j})=y^j(-\tau_rH^{\ r}_{i\
 jk}-\tau_{i.r}R^r_{\ jk})=-\tau_rH^{\ r}_{i\
 0k}-\tau_{i.r}R^r_{\ k}.
  \ee
Now by $R^i_{\ k}=R\delta^i_k-\tau_ky^i$ we have
 \be\label{yr65}
 H^{\ r}_{i\
jk}=\frac{1}{3}\big[R_{.j.i}\delta^r_k-\tau_{k.j.i}y^r-\tau_{k.j}\delta^r_i-\tau_{k.i}\delta^r_j
-(j/k)\big].
 \ee
 Plugging (\ref{yr65}), $R^i_{\ k}=R\delta^i_k-\tau_ky^i$,
 $\tau_{i;j}=0$ and $R=0$ into (\ref{yr64}) we obtain
 $\tau_i\tau_k=0$, which gives $\tau_i=0$. This gives the proof of
 (\ref{ycw65}).

 Assume  $R_{;i}=0 (R\ne 0)$. By (\ref{ycw64}) we have
 $R_{.i}=2\tau_i$. So
 $\tau_{i;k}=\frac{1}{2}R_{.i;k}=\frac{1}{2}R_{;k.i}=0$. This gives the proof of
 (\ref{ycw66}).  \qed

\

{\it Proof of Theorem \ref{th01} :} (i) Assume that {\bf G} is of
constant curvature. Then we have $R^i_{\ k}=R\delta^i_k-\tau_ky^i$
with $\tau_{i;k}=0$ by definition. By (\ref{ycw65}) we immediately
obtain the desired conclusion since $Ric=(n-1)R$. Conversely, let
the Riemann curvature be zero or the Ricci curvature $Ric$ satisfy
$Ric_{;i}=0 (Ric\ne0)$. If $R^i_{\ k}=0$, then $\tau_i=0$ and so
$\tau_{i;k}=0$. If $Ric_{;i}=0 (Ric\ne0)$, then we have $R_{;i}=0\
(R\ne0)$. So by (\ref{ycw64}) we have $R_{.i}=2\tau_i$. Thus we
obtain $\tau_{i;k}=\frac{1}{2}R_{.i;k}=\frac{1}{2}R_{;k.i}=0$. By
definition, {\bf G} is of constant curvature.

(ii) It follows directly from (\ref{ycw65}) and  (\ref{ycw64}),
and the definitions for  a spray of isotropic curvature and
constant curvature.

(iii) If a Finsler metric $L$ is of constant flag curvature
$\lambda$, then its spray has the Riemann curvature $R^i_{\
k}=\lambda (L\delta^i_k-y_ky^i)$. So its spray is of scalar
curvature $R^i_{\ k}=R\delta^i_k-\tau_ky^i$ with $R=\lambda L$ and
$\tau_k=\lambda y_k$. Thus we have $\tau_{i;k}=0$, which implies
that the spray is of constant curvature. Conversely, if the spray
{\bf G} of a Finsler metric $L$ is of constant curvature, then the
Riemann curvature of {\bf G} has the form $R^i_{\
k}=R\delta^i_k-\tau_ky^i$ with $\tau_{i;k}=0$. So $L$ is of scalar
flag curvature (put the flag curvature as $\lambda$). We have
$R=\lambda L$ and $\tau_k=\lambda y_k$. Thus we get
$\lambda_{;k}=0$ by $\tau_{i;k}=0$. If $\lambda=0$, then $L$ has
constant flag curvature $0$. If $\lambda\ne 0$, then since
 $\lambda$ is a homogeneous function of degree
$0$ satisfying $\lambda_{;k}=0$,  we immediately obtain
$\lambda=constant$ by putting $T=\lambda, p=0$ in (\ref{yg63}) of
Proposition \ref{prop32}.

(iv) See Examples \ref{ex73} and \ref{ex74} below. \qed

\

{\it Proof of Theorem \ref{th001} :} The spray {\bf G} has the
Riemann curvature $R^i_{\ k}=R\delta^i_k-\tau_ky^i$ with
$\tau_{i;k}=0$ (note that $Ric=(n-1)R$). By Theorem
\ref{th01}(ii), the spray {\bf G} is of isotropic curvature since
{\bf G} is of constant curvature.

If {\bf G} is (locally) Finsler-metrizable, then by Theorem
\ref{th1}, it is clear that there has $Ric=0$ or $Ric$ is a
Finsler metric when $Ric\ne 0$.

Conversely, if $Ric=0$, then by Theorem \ref{th1}(i), {\bf G} is
(locally) Finsler-metrizable. If $Ric\ne0$, then by assumption
$Ric$ is a Finsler metric. It follows from Theorem \ref{th01}(i)
that $Ric_{;i}=0$. By Theorem \ref{th1}(iii), {\bf G} is (locally)
Finsler-metrizable.   \qed

\section{Locally projectively flat sprays}

In this section, we consider the properties and metrizability of
locally projectively flat sprays on a manifold. A locally
projectively flat spray is always of scalar curvature. Locally, we
let {\bf G} be a projectively flat spray with $G^i=Py^i$ defined
on an open set of $R^n$. Then by (\ref{y004}), the Riemann
curvature tensor $R^i_{\ k}$ is in the form $R^i_{\
k}=R\delta^i_k-\tau_ky^i$ with
 \be\label{ygj31}
 R=P^2-P_{x^r}y^r,\ \ \tau_k=PP_{y^k}+P_{x^ry^k}y^r-2P_{x^k}.
 \ee

 \subsection{Some basic results}

\begin{lem}\label{lem51}
Let $G^i=Py^i$ be a spray defined on an open set of $R^n$. Then
 (\ref{ygj31}) becomes
 \be\label{ygj32}
 R=-P^2-P_{;0},\ \ \ \ \tau_k=P_{k;0}-2P_{;k}-PP_k, \ \ \
 (P_k:=P_{.k}).
 \ee
\end{lem}

\begin{lem}\label{lem52}
 In (\ref{ygj32}), $R=0$ and $\tau_k=0$  $\Longleftrightarrow$
 $P_{;k}+PP_k=0$. So the spray $G^i=Py^i$ has vanishing Riemann
 curvature iff. $P_{;k}+PP_k=0$.
 \end{lem}

 {\it Proof :} If $R=0$ and $\tau_k=0$, then $P_{;0}=-P^2$ implies
 $P_{;k}+P_{k;0}=-2PP_k$. Further by $P_{k;0}=2P_{;k}+PP_k$ we
 obtain $P_{;k}+PP_k=0$. Conversely, let $P_{;k}+PP_k=0$. Then we
 have $P_{;0}=-P^2$ $(R=0)$. Now again $P_{;0}=-P^2$ implies
 $P_{;k}+P_{k;0}=-2PP_k$, which implies $\tau_k=0$ by
 $P_{;k}+PP_k=0$. \qed

 \

By Theorem \ref{th1}(i) and Lemma \ref{lem52}, a spray $G^i=Py^i$
satisfying $P_{;k}+PP_k=0$ is locally Finsler-metrizable (also see
\cite{Shen5, LS}). By (\ref{ygj32}) we can easily obtain the
following lemma.
 \begin{lem}\label{lem53}
Let $G^i=Py^i$ be a spray with $R^i_{\ k}=R\delta^i_k-\tau_ky^i$.
Then we have
 \be\label{ygj33}
 R_{.i}-2\tau_i=3(P_{;i}-P_{i;0}).
 \ee
 So {\bf G} is of isotropic curvature iff. $P_{;i}=P_{i;0}$.
 \end{lem}

By Lemma \ref{lem53}, we can easily obtain the following
Proposition \ref{prop64}.

\begin{prop}\label{prop64}
 Let $G^i=Py^i$ be a Berwald spray. Then {\bf G} is of isotropic
 curvature iff. $P$ is a (local) exact 1-form given by
 $P=\sigma_{x^i}y^i$ for a scalar function $\sigma=\sigma(x)$.
\end{prop}

The following Proposition \ref{prop65} directly follows from
 Lemma \ref{lem53} and Theorem \ref{th1}.

\begin{prop}\label{prop65}
Let $G^i=Py^i$ ($P\ne 0$) be a spray defined on an open set of
$R^n$ with $P_{;i}=0$. Then {\bf G} is of
 isotropic curvature. Further {\bf G} is
     Finsler metrizable iff. $P$ is a Finsler metric, and in this
     case, $P^2$ is a Finsler metric of constant flag curvature
     $-1$.
\end{prop}

 By Theorem \ref{th3}(ia) and (\ref{ygj33}) we have the following
Proposition \ref{prop66} (cf. \cite{LMY}).

\begin{prop}\label{prop66}
 For an $n(\ge3)$-dimensional Berwald spray {\bf G} with $G^i=Py^i$, if {\bf G} is not of
 isotropic curvature (or equivalently, $P$ is not a closed 1-form), then {\bf G} is not Finsler-metrizable.
\end{prop}

\subsection{Proof of Theorem \ref{th2}}
 Now we give the proof of Theorem \ref{th2} as follows.

 If  the spray {\bf G} is given
by (\ref{ygj1}) (defined on $\mathcal{C}(U)$), then {\bf G} is
Finsler-metrizable (on $\mathcal{C}(U)$) since it is easy to check
that the  metric $L$ given by (\ref{ygj2}) induces the spray {\bf
G}. Further, a direct computation shows that the Ricci curvature
$Ric$ of {\bf G} is equal to $(n-1)L^*$ (on $\mathcal{C}(U)$)
since {\bf G} is given by (\ref{ygj1}). So $L:=Ric/(n-1)$ is a
Riemann metric  and $L$ induces {\bf G}.

Conversely, suppose that the spray {\bf G} is Finsler-metrizable
(on $\mathcal{C}(U)$). Since {\bf G} is a projectively flat
Berwald spray of isotropic curvature, it follows from Proposition
\ref{prop64} that  {\bf G} has the form $G^i=\sigma_0y^i$, where
$\sigma_i=\sigma_{x^i}$ is the differential of some scalar
function $\sigma=\sigma(x)$. Since $Ric\ne 0$ and {\bf G} is
Finsler-metrizable, $Ric$ is a Finsler metric and $Ric_{;i}=0$ by
theorem \ref{th1} and the fact  that
 a locally  projectively flat Finsler
metric of isotropic flag curvature is of constant flag curvature
(\cite{Ber0}). Since $Ric=(n-1)[-(\sigma_0)^2-\sigma_{0;0}]$ (see
(\ref{ygj32})) and $Ric_{;i}=0$, we have
 \be\label{ygj34}
 \big[(\sigma_0)^2+\sigma_{0;0}\big]_{;i}=0.
 \ee
Now we only need to solve the scalar function $\sigma$ from
(\ref{ygj34}). It is easy to see that (\ref{ygj34}) is equivalent
to
 \be\label{ygj35}
(\sigma_j\sigma_k+\sigma_{j;k})_{;i}=0.
 \ee
Now by $G^i=\sigma_0y^i$ we have
 $$
 G^i_j=\sigma_jy^i+\sigma_0\delta^i_j,\ \ \ \  G^i_{jk}=\sigma_j\delta^i_k+\sigma_k\delta^i_j.
 $$
For convenience, we put
$\sigma_{ij}:=\sigma_{x^ix^j},\sigma_{ijk}:=\sigma_{x^ix^jx^k}$.
By a direct computation we obtain
 \beq
\sigma_{i;j}&&\hspace{-0.6cm}=\sigma_{ij}-2\sigma_i\sigma_j,\ \ \
\ \  \sigma_i\sigma_j+\sigma_{i;j}
=\sigma_{ij}-\sigma_i\sigma_j,\nonumber \\
 (\sigma_j\sigma_k+\sigma_{j;k})_{;i}&&\hspace{-0.6cm}
 =\sigma_{ijk}-2(\sigma_i\sigma_{jk}+\sigma_j\sigma_{ik}+\sigma_k\sigma_{ij})
 +4\sigma_i\sigma_j\sigma_k.\label{ygj36}
 \eeq
Putting $u=e^{-2\sigma}$, it follows from (\ref{ygj36}) that
(\ref{ygj35}) is equivalent to $u_{ijk}=0$. So $u$ is a polynomial
in ($x^i)$ of degree two. Thus the scalar function $\sigma$
satisfying (\ref{ygj35}) is given by
 \be\label{ycw51}
 \sigma(x)=-\frac{1}{2}\ln|x'Ax+\langle B,x\rangle+C|,
 \ee
 where $'$ denotes the transpose of a matrix, $A$ is a constant symmetric matrix, $B$ is a constant vector
  and $C$ is a constant number. By $G^i=\sigma_0y^i$ and (\ref{ycw51}), we see that {\bf G} is locally given by (\ref{ygj1}).
   Meanwhile, by (\ref{ycw51}), we have
 \be\label{ygj37}
  -(\sigma_0)^2-\sigma_{0;0}=\frac{4(x'Ax+\langle
  B,x\rangle+C)y'Ay-(2x'Ay+\langle
  B,y\rangle)^2}{4(x'Ax+\langle
  B,x\rangle+C)^2}.
  \ee
Since $R(=-(\sigma_0)^2-\sigma_{0;0})\ne 0$ by assumption, by
Theorem \ref{th1}, $R$ is a metric, which means that $A,B,C$
should satisfy certain condition such that $L$ given by
(\ref{ygj2}) is a metric.

 A simple
computation shows that if $det(A)\ne 0$, then $L$ given by
(\ref{ygj2}) is a metric iff. the constant quantities $A,B,C$
satisfy $4C\ne B'A^{-1}B$.    \qed

\begin{rem}\label{rem57}
  If $-(\sigma_0)^2-\sigma_{0;0}=0$, we can obtain $\sigma=-\ln|\langle
B,x\rangle+C|$. In this case,   {\bf G} has vanishing Riemann
curvature, and it can be locally induced by a Finsler metric $L$.
It is clear that $L$ is a Minkowskin metric.
\end{rem}

The above proof and Remark \ref{rem57} have actually given the
proof of the following theorem for the local structure of a
locally projectively flat Berwald spray with constant curvature.

\begin{thm}\label{th68}
 Let {\bf G} be a locally projectively flat Berwald spray
 on a manifold $M$. Then {\bf G} is of constant curvature if and only if  {\bf G}
can be locally expressed as (\ref{ygj1}).
\end{thm}

\section{Examples}\label{sec7}

 As an application of Theorem \ref{th1},
Theorem \ref{th2} gives the necessary and sufficient condition for
the Berwald spray $G^i=(\sigma_{x^r}y^r)y^i$ to be
Finsler-metrizable. In this section, we are going to give more
examples to support some of the basic theories we introduce in the
above sections.

\begin{ex}
 Let ${\bf G}$ be a two-dimensional Berwald spray given by
  $$
 G^1=f(x^1)[(y^1)^2-(y^2)^2], \ \ \ \ \ G^2=2f(x^1)y^1y^2,
  $$
  where $f(x^1)$ is a non-constant function. It can be directly
  verified that ${\bf G}$ is of isotropic curvature  $R^i_{\
  k}=R\delta^i_k-\frac{1}{2}R_{.k}y^i$ satisfying
   $$R=-2f'(x^1)[(y^1)^2+(y^2)^2],\ \ \ \ R_{;i}=R\omega_i$$
   $$
 \omega_1:=\frac{f''(x^1)}{f'(x^1)}-4f(x^1)=\big[\ln|f'(x^1)e^{-4\int
 f(x^1)dx^1}|\big]_{x^1}=[\ln|\lambda|]_{x^1},\ \ \ \ \ \omega_2:=0.
   $$
   Then by Theorem \ref{th1}, ${\bf G}$ is Finsler-metrizable induced by
   the Riemann metric
   $$L:=R/\lambda=-2e^{4\int
 f(x^1)dx^1}[(y^1)^2+(y^2)^2]
 $$ of
   isotropic sectional curvature $$\lambda=f'(x^1)e^{-4\int
 f(x^1)dx^1}.$$
\end{ex}

\begin{ex}
Let ${\bf G}$ be a two-dimensional  spray given by
  $$
 G^1=\frac{1}{2r}y^2\sqrt{(y^1)^2+(y^2)^2}, \ \ \ \ \ G^2=-\frac{1}{2r}y^1\sqrt{(y^1)^2+(y^2)^2},
  $$
 where $r$ is a constant. This spray appears in Example 4.1.3 of \cite{Shen6}. From the
  completeness of geodesics of ${\bf G}$ and Hopf-Rinow theorem, it is concluded that ${\bf
  G}$ is not globally Finsler-metrizable. But actually, it is
  even not locally Finsler-metrizable by Theorem \ref{th1} (Example 12.4.1 in \cite{Shen6}
   shows it is locally projectively  Finsler-metrizable).

  By a direct computation, ${\bf G}$ is of isotropic curvature  $R^i_{\
  k}=R\delta^i_k-\frac{1}{2}R_{.k}y^i$ satisfying
   $$
 R=r^{-2}\big[(y^1)^2+(y^2)^2\big],\ \ \ \ R_{;1}=r^{-3}y^2\sqrt{(y^1)^2+(y^2)^2},\ \ \ \
 R_{;2}=-r^{-3}y^1\sqrt{(y^1)^2+(y^2)^2}.
   $$
   Therefore, by Theorem \ref{th1}, ${\bf G}$  is  not locally
   Finsler-metrizable anywhere since there exists no closed 1-form
   $\omega=\omega_i(x)dx^i$ such that $R_{;i}=R\omega_i$.
\end{ex}

\begin{ex}\label{ex73}
 Let ${\bf G}_c$ be a Berwald spray given by
  $$
 G^i_c=\frac{-|y|^2x^i+c\langle x,y\rangle y^i}{1-|x|^2},
  $$
  where $c$ is a constant.
The spray ${\bf G}_2$ is
  induced by the  Riemann metric
  $L=4|y|^2/(1-|x|^2)^2$ of constant sectional curvature $-1$, and by computation, we can see that ${\bf G}_1$
  is just the spray given by Example 4.1.4 in \cite{Shen6}. ${\bf G}_c$
  actually is  projectively flat
  although at present coordinate it is not in the form $Py^i$.

 By a direct computation,  ${\bf G}_c$ is of isotropic curvature  $R^i_{\ k}=R\delta^i_k-\frac{1}{2}R_{.k}y^i$
 with
  $$
R=-\frac{\big[(c-2)|x|^2+c+2\big]|y|^2-c(c-2)\langle
x,y\rangle^2}{(1-|x|^2)^2}.
  $$
  Further, we have
   $$
 R_{;i}=\frac{2(c-2)}{(1-|x|^2)^3}\Big\{\big[(c-1)|x|^2+c+1\big](|y|^2x^i+2\langle
 x,y\rangle y^i)-2c(c-1)\langle
 x,y\rangle^2x^i\Big\}.
   $$

   (i) By Theorem \ref{th1}, it is easy to see that ${\bf G}_c$ is
  Finsler-metrizable iff. $c=2$. When $c=2$, we obtain the Riemann
  metric $L:=R=-4|y|^2/(1-|x|^2)^2$ of constant
  curvature 1.

  (ii) ${\bf G}_c$ is of constant  curvature iff. $c=2$.

  (iii) When $c\ne2$, ${\bf G}_c$ is of isotropic  curvature but not
of constant curvature in any dimension. This fact is different
from the Finslerian case.
\end{ex}

\begin{ex}\label{ex74}
Let ${\bf G}$ be a spray given by $G^i=\sigma_0y^i$ with
$\sigma_i:=\sigma_{x^i}, \sigma=\sigma(x)$.

(i) ${\bf G}$ is of isotropic curvature for any $\sigma$ but not
of constant curvature by Theorem \ref{th68} if
   $$
 \sigma(x)\ne -\frac{1}{2}\ln|x'Ax+\langle B,x\rangle+C|.
 $$

 (ii) Taking $A=(\delta_{ij}),B=0,C=0$ in (\ref{ygj1}) and
 (\ref{ygj2}), we have
  \be\label{cw52}
 G^i=-\frac{\langle x,y\rangle}{|x|^2}y^i,\ \ \ \
 L=\frac{|x|^2|y|^2-\langle x,y\rangle^2}{|x|^4}.
  \ee
  In (\ref{cw52}), $L$ is not a metric since $(L_{.i.j})$ is
  degenerate. By Theorem \ref{th2}, or Theorem \ref{th1}(ii)
  ($R=L$ is not a Finsler metric), the spray {\bf G} given by
  (\ref{cw52}) is not Finsler-metrizable, although this spray {\bf
  G} satisfies (\ref{cw16}).
\end{ex}

\begin{ex}
 Let {\bf G} be a spray on $R^3$ given by (\cite{Shen6}, Example 4.1.2)
 $$
G^1=0,\ \ G^2=x^1(y^1)^2+x^3(y^3)^2,\ \ G^3=0.
 $$
 {\bf G} has zero Riemann curvature. Then by Theorem \ref{th1}(i),
 {\bf G} is locally Finsler-metrizable. Actually, if the following
 function is a Finsler metric
  $$
 L(x,y)=\Big[f\Big(\frac{y^3}{y^1},(x^1)^2+\frac{y^2}{y^1}+(x^3)^2\frac{y^3}{y^1}\Big)y^1\Big]^2,
  $$
  then $L$ induces {\bf G} on some open set of $R^3$, and $F$ should be a locally
  Minkowski metric since {\bf G} is a Berwald spray (with zero Riemann
  curvature). Consider a special case
  with $f(r,s)=\sqrt{r^2+s}$ and we obtain
   $$L(x,y)=(y^3)^2+(x^1)^2(y^1)^2+y^1y^2+(x^3)^2y^1y^3,$$
   which is a singular Riemann metric globally defined on $R^3$. So the spray {\bf G}
   is globally Finsler-metrizable.
\end{ex}

\begin{ex}\label{ex76}
 Let {\bf G} be a two-dimensional spray given by
  $$
 G^1=-f(x^1)g'(t)(y^1)^2,\ \ \  G^2=f(x^1)\big[g(t)y^1-g'(t)y^2\big]y^1,\ \ \ \ (t:=y^2/y^1).
  $$
 where $f, g$ are two smooth functions. A direct computation gives $R^i_{\
  k}=R\delta^i_k-\tau_ky^i$ with
   $$
   R=0,\ \ \  \tau_1=Ay^2,\ \ \
 \tau_2=-Ay^1,\ \ \  A:=2f^2(x^1)g(t)g'''(t)-f'(x^1)g''(t).
   $$
  We can choose the two functions $f, g$ satisfying $A\ne 0$. So in this case, {\bf G} is not Finsler
   metrizable by Theorem \ref{th32}(ii).
\end{ex}

\vspace{0.5cm}

\noindent Guojun Yang \\
Department of Mathematics \\
Sichuan University \\
Chengdu 610064, P. R. China \\
yangguojun@scu.edu.cn

\end{document}